\documentclass{amsart}

\usepackage{amsmath,amsthm,amssymb}

\newtheorem{thm}{Theorem}[section]
\newtheorem{id}[thm]{Identity}

 \theoremstyle{definition}
 
 \theoremstyle{remark}
 \newtheorem*{rem}{Remark}
 
 \numberwithin{equation}{section}

\newcommand{\qbin}[3]{\genfrac{[}{]}{0pt}{}{#1}{#2}_{#3}}
\newcommand{\agp}[4]{\left[ \genfrac{}{}{0pt}{}{#1}{#2} ; #3, #4 \right]}

\begin{document}
\title[Two-variable Ramanujan-type identities]{Polynomial Generalizations of two-variable Ramanujan type identities}

\author{James McLaughlin}
\email{jmclaughlin@wcupa.edu}
\address{Department of Mathematics, West Chester University, West Chester, Pennsylvania }

\author{Andrew V. Sills}
\email{ASills@GeorgiaSouthern.edu}
\address{Department of Mathematical Sciences, Georgia Southern University, Statesboro,
Georgia
30460-8093}

\date{\today}

\dedicatory{Dedicated to Doron Zeilberger on the occasion of his sixtieth birthday.}


\begin{abstract}
We provide finite analogs of a pair of two-variable $q$-series
identities from Ramanujan's lost notebook and a companion identity.
\end{abstract}

\maketitle

``The progress of mathematics can be viewed as progress from the infinite to the finite."
---Gian-Carlo Rota (1983)

\section{Introduction}
At the top of a page in the lost notebook~\cite[p. 33]{R88} (cf.~\cite[p. 99, Entry 5.3.1]{AB09}):, Ramanujan recorded an identity equivalent to the following:
\begin{equation} \label{R1}
  \sum_{j=0}^\infty \frac{q^{2j^2} (zq; q^2)_j (q/z;q^2)_j}{(q^2;q^2)_{2j}}
 = \frac{(zq^3, q^3/z, q^6;q^6)_\infty}{(q^2;q^2)_\infty},
\end{equation}
where we are employing the standard notation for rising $q$-factorials,
\[ (A;q)_\infty := (1-A)(1-Aq)(1-Aq^2)\cdots
 \mbox{ and }(A;q)_n:=  \frac{(A;q)_\infty}{(Aq^n;q)_\infty},\] and
 \[ (A_1, A_2, \cdots, A_r; q)_\infty := (A_1;q)_\infty (A_2;q)_\infty \cdots (A_r;q)_\infty. \]

 In a recent paper~\cite{MSZ09}, we found a partner to~\eqref{R1} that
Ramanujan appears to have missed:
\begin{equation} \label{R1partner}
\sum_{j=0}^\infty \frac{ q^{j(j+1)} (z;q)_j (q/z;q)_{j+1}}{(q;q)_{2j+1}}
= \frac{(zq^2, q/z, q^3; q^3)_\infty}{(q;q)_\infty}.
\end{equation}

Later on the same page of the lost notebook, Ramanujan recorded
\cite[p. 103, Entry 5.3.5]{AB09}
\begin{equation}\label{R2}
 \sum_{j=0}^\infty \frac{q^{j^2} (zq;q^2)_j (q/z;q^2)_j }{(q;q^2)_j (q^4;q^4)_j}
 = \frac{(zq^2, q^2/z, q^4;q^4)_\infty (-q;q^2)_\infty}{(q^2;q^2)_\infty}.
\end{equation}

For further discussion of these three identities, see~\cite{MS10}.

\begin{rem}
Out of respect for Doron's ultra-finitist philosophy, we deliberately refrain from stating conditions
on $q$ and $z$ which imply analytic convergence of the infinite series and products in
~\eqref{R1}--\eqref{R2}.
\end{rem}

 The preceding identities stand out among identities of Rogers-Ramanujan type
because they are \emph{two-variable} series-product identities.
While Rogers-Ramanujan type identities admit two-variable generalizations,
most lose the infinite product representation in the two-variable case.

 For example, in the standard two
variable generalization of the first Rogers-Ramanujan identity,
\begin{equation} \label{aRR1}
\sum_{j=0}^\infty \frac{z^j q^{j^2}}{(q;q)_j}
= \frac{1}{(zq;q)_\infty} \sum_{n=0}^\infty \frac{(-1)^j z^{2j}
q^{j(5j-1)/2} (1-zq^{2j}) (z;q)_j}{ (1-z)(q;q)_j},
\end{equation}
the right hand side reduces to an infinite product \emph{only} for certain
particular values of $z$, e.g.
$z=1$ gives the first Rogers-Ramanujan identity,
\begin{equation}
  \sum_{j=0}^\infty \frac{q^{j^2}}{(q;q)_j} = \frac{1}{(q;q^5)_\infty (q^4;q^5)_\infty},
\end{equation} while $z=q$ gives the second Rogers-Ramanujan identity,
\begin{equation}
  \sum_{j=0}^\infty \frac{q^{j(j+1)}}{(q;q)_j} = \frac{1}{(q^2;q^5)_\infty (q^3;q^5)_\infty},
\end{equation}
after application of the Jacobi triple product identity~\cite[p. 17, Eq. (1.4.8)]{AB09}.

In~\cite[\S3]{S03}, the second author presented nontrivial polynomial generalizations
of all 130 Rogers-Ramanujan type identities appearing in Slater's paper~\cite{S52}.
All of Slater's identities involved one variable only.
Here, we demonstrate that the methods employed in~\cite{S03} can be used to obtain
polynomial generalizations of the rarer species of two-variable $q$-series-product
identities as well.

\section{Polynomial Generalizations}
Define the standard binomial co\"efficient by
\[  \qbin{A}{B}{q} := \left\{
 \begin{array}{ll}
  \displaystyle{\frac{(q;q)_A}{(q;q)_B (q;q)_{A-B}}}, &\mbox{ if $0\leq B \leq A$} \\
  0, &\mbox{otherwise}
 \end{array}, \right. \]
and the modified $q$-binomial co\"efficient by
\[  \qbin{A}{B}{q}^* := \left\{
 \begin{array}{ll}
  1, &\mbox{ if $A=-1$ and $B=0$,} \\
  \qbin{A}{B}{q}, &\mbox{otherwise}
 \end{array}. \right. \]

In~\cite{AB87}, Andrews and Baxter define several $q$-analogs of trinomial
co\"efficients; we shall require one
 of them here:
\[ T_0 (L, A; q) := \sum_{r=0}^L (-1)^r \qbin{L}{r}{q^2} \qbin{2L-2r}{L-A-r}{q}. \]

  More recently, Andrews~\cite{A05} introduced the following generalization of
the $q$-binomial co\"efficient:
\[
\agp{A}{B}{q}{z} := \left\{
\begin{array}{ll}
   0 & \mbox{if $B<0$}\\
   1 & \mbox{if $B=0$ or $B=A$}\\
   \sum_{h=0}^B z^h \qbin{A-B+h-1}{h}{q} & \mbox{if $0<B<A$} \\
   (zq^{A-B};q)_{B-A} & \mbox{if $B>A$.}
   \end{array}
\right.
\]

The following polynomial generalizations of~\eqref{aRR1} are known:
\begin{equation} \label{aRRf1}
\sum_{j=0}^n z^j q^{j^2} \qbin{n}{j}{q} =
\sum_{j=0}^n (-1)^j q^{j(5j-1)/2} (1-z q^{2j}) \qbin{n}{j}{q} \frac{1}{(zq^j;q)_{n+1}}
\end{equation}  (see~\cite{A74,B81a,ET90,P94}),
\begin{multline} \label{aRRf2}
\sum_{j=0}^n z^n q^{n^2} = \sum_{0\leq 2j \leq n} (-1)^j z^{2j} q^{j(5j-1)/2} (1-zq^{2j})
 \qbin{n}{j}{q} \qbin{n-j}{j}{q} (q;q)_j \\ \times
 \frac{(z^2 q^{n+2j+1};q)_{n-2j}}{(zq^j;q)_{n-j+1}}
\end{multline} ~\cite[Eq. (3.5)]{B81b}, and
\begin{multline} \label{aRRf3}
\sum_{j=0}^n z^j q^{j^2} \agp{n}{j}{q}{q} \\
= \sum_{0\leq 2j \leq n} (-1)^j z^{2j} q^{j(5j-1)/2} \agp{n}{j}{q}{q} \agp{2n+1-2j}{n-2j}{q}{zq^j}\\
-\sum_{0\leq 2j \leq n-1} (-1)^j z^{2j+1} q^{j(5j+3)/2} \agp njqq \agp{2n-2j}{n-2j-1}{q}{zq^j},
\end{multline}
\cite[p. 41, Eq. (1.11)]{A05}.

Andrews~\cite{A05} notes that one of his motivations for introducing~\eqref{aRRf3} is
that both sides of the equation are clearly polynomials term by term, whereas this is
\emph{not} the case for the right hand sides of
\eqref{aRRf1} and~\eqref{aRRf2}.  The polynomial identities we introduce below
also have this desirable feature.

Notice that in each of the identities below, the summands have finite support, and follow
the natural bounds (i.e. each summation could be taken over all integers, and no
nonzero terms would be added).

\begin{id}[Polynomial Generalization of~\eqref{R1}] \label{R1finite}
\begin{multline}
\sum_{j=0}^{\lfloor n/2 \rfloor} \sum_{h=0}^j \sum_{i=0}^j  (-1)^{h+i} z^{h-i} q^{h^2+i^2+2j^2}
\qbin{j}{h}{q^2} \qbin{j}{i}{q^2} \qbin{j+ \lfloor \frac{n-h-i}{2} \rfloor}{2j}{q^2}\\
 =\sum_{j=-\infty}^\infty (-1)^j z^j q^{3j^2} \qbin{n-1}{\lfloor \frac{n+3j-1}{2} \rfloor}{q^2}
 +\epsilon_n (z,q),
 \end{multline}
 where
 \begin{equation}
 \epsilon_n (z,q) \\= \left\{
 \begin{array}{ll}
\displaystyle{ \sum_{j=-\infty}^\infty z^{2j} q^{12j^2+6j+n} \qbin{n-1}{\frac{n+6j}{2}}{q^2}^*}  &\mbox{if $2\mid n$} \\
\displaystyle{- \sum_{j=-\infty}^\infty z^{2j-1} q^{12j^2-6j+n} \qbin{n-1}{\frac{n+6j-3}{2}}{q^2}^*} &\mbox{if $2\nmid n$}
\end{array} \right.
 \end{equation}.
\end{id}

\begin{id}[Polynomial Generalization of~\eqref{R1partner}]\label{R1PartnerFinite}
\begin{multline}
\sum_{j=0}^{\lfloor n/2 \rfloor} \ \sum_{h=0}^j \sum_{i=0}^{j+1} (-1)^{h+i} z^{h-i} q^{\binom{h}{2} +\binom{i+1}{2} +j(j+1) }
\qbin{j}{h}{q} \qbin{j+1}{i}{q}
\qbin{j+1 + \lfloor \frac{n-h-i}{2} \rfloor }{2j+1}{q} \\
= \sum_{j=-\infty}^\infty (-1)^j z^j q^{j(3j+1)/2} \qbin{n}{ \lfloor \frac{n+3j+2}{2}  \rfloor}{q}
+ \epsilon_n(z,q),
\end{multline}
where
\begin{equation}
\epsilon_n(z,q) = \left\{
  \begin{array}{ll}
    \displaystyle{ \sum_{j=-\infty}^\infty z^{2j} q^{6j^2 - 2j + n/2} \qbin{n}{\frac n2 + 3j}{q} }&\mbox{if $2\mid n$},\\
     \displaystyle{-\sum_{j=-\infty}^\infty z^{2j+1} q^{6j^2 + 4j +\frac 12 + \frac n2} \qbin{n}{\frac{n+6j+3}{2}}{q}}
       &\mbox{if $2\nmid n$}.
  \end{array}
\right.
\end{equation}
\end{id}

\begin{id}[Polynomial Generalization of~\eqref{R2}] \label{R2finite}
\begin{multline}
 \sum_{j=0}^{n} \ \sum_{h=0}^j \sum_{i=0}^j \sum_{\ell=0}^{n-h-i}
  (-1)^{h+i+\ell} z^{h-i}q^{h^2+i^2+j^2+2\ell}
  \qbin jh{q^2} \qbin ji{q^2} \qbin{j+\ell-1}{\ell}{q^2}^*
  \\ \times
\qbin{n-h-i+j-\ell}{2j}{q}
  \\ = \sum_{j=-\infty}^{\infty} (-1)^j z^j q^{2j^2} \bigg( T_0(n, 2j; q) + T_0 (n-1, 2j; q) \bigg)
\end{multline}
\end{id}

\section{Derivation and a method of proof}
\subsection{Identity~\ref{R1finite}}
Recall the following consequences of the $q$-binomial theorem:
\begin{equation} \label{qBT1}
  (t;q)_j = \sum_{h=0}^j (-1)^h t^h  q^{h(h-1)/2} \qbin{j}{h}{q}
\end{equation}
\begin{equation} \label{qBT2}
  \frac{1}{(t;q)_j} = \sum_{h=0}^\infty t^h \qbin{h+j-1}{h}{q}^*
\end{equation}

The derivation of~\eqref{R1finite} is via the method used for
the derivations of polynomial versions of Rogers-Ramanujan
type identities (in $q$ only) as introduced by Andrews~\cite[Chapter 9]{A86},
and further explored by Santos~\cite{S91} and the second author~\cite{S03, S04}.
We shall consider the details of~\eqref{R1} only; ~\eqref{R1partner} and~\eqref{R2}
may be treated analogously.

We begin with the left hand side of~\eqref{R1}
\begin{equation} \label{R1lhs}
 \phi(z,q):=\sum_{j=0}^\infty \frac{q^{2j^2} (zq;q^2)_j (q/z;q^2)_j}{(q^2;q^2)_{2j} }.
 \end{equation}
Now define the following generalization of $\phi(z,q)$:
\begin{equation}
  f(t):= f(t; z, q) := \sum_{j=0}^\infty \frac{t^{2j} (1+t) q^{2j^2} (tzq;q^2)_j (tq/z;q^2)_j}
  { (t^2;q^2)_{2j+1}},
  \end{equation}
and let $P_n (z,q)$ be defined by
\[ f(t) = \sum_{n=0}^\infty P_n(z,q) t^n . \]
Note that
\[ \lim_{t\to 1-} (1-t)f(t; z,q) = \phi (z,q) \] and
\[ \lim_{n\to\infty} P_n (z,q) = \phi(z,q). \]

\begin{align*}
f(t) &= \sum_{j=0}^\infty \frac{t^{2j} (1+t) q^{2j^2} (tzq;q^2)_j (tq/z;q^2)_j}
  {(t^2;q^2)_{2j+1}}\\
  & = \frac{1}{1-t}+\sum_{j=1}^\infty \frac{t^{2j} (1+t) q^{2j^2} (tzq;q^2)_j (tq/z;q^2)_j}
  {(t^2;q^2)_{2j+1}}\\
 & =  \frac{1}{1-t}+\sum_{j=0}^\infty \frac{t^{2j+2} q^{2j^2+4j+2} (tzq;q^2)_{j+1}
 (tq/z;q^2)_{j+1}}
  { (t^2;q^2)_{2j+3}}\\
  & = \frac{1}{1-t} + \frac{t^2 q^2 (1-tzq)(1-tq/z)}{(1-t^2 q^2)(1+tq^2)(1-t)} f(tq^2).
\end{align*}
Thus,
\begin{equation*}
(1-t^2q^2)(1+tq^2)(1-t) f(t) = (1-t^2q^2)(1+tq^2) + t^2 q^2(1-tzq)(1-tq/z) f(tq^2),
\end{equation*}
which immediately implies
\begin{multline} \label{qDE}
f(t) = (1+ tq^2 - t^2 q^2 - t^3q^4) + \left( (1-q^2)t + 2q^2 t^2 + (q^4-q^2)t^3 - q^4t^4
\right)f(t) \\
+ \left( q^2 t^2 - (z+z^{-1})q^3 t^3 + q^4 t^4 \right) f(t q^2).
\end{multline}
Upon recalling that $f(t) = \sum_{n=0}^\infty P_n(z,q) t^n$, and extracting the
co\"efficients of $t^n$ from~\eqref{qDE},
we find that the  $P_n = P_n(z,q)$ satisfy the fourth order recurrence
\begin{multline} \label{R1rec}
P_n = (1-q^2) P_{n-1} + (2q^2+q^{2n-2}) P_{n-2} +
\left( q^4 - q^2- (z+z^{-1})q^{2n-3}\right)P_{n-3} \\+ (q^{2n-4} - q^4) P_{n-4}
\end{multline}
with initial conditions
\begin{equation} \label{R1ic}
P_0 = P_1 = 1; \qquad P_2 = 1+q^2; \qquad P_3 = 1+q^2 - (z+z^{-1}) q^3.
\end{equation}
Thus we now have a full characterization of the $P_n(z,q)$ via a recurrence with
initial conditions.

   Next, we use $f(t)$ to derive the left hand side of~\eqref{R1finite}.
{\allowdisplaybreaks
\begin{align*}
\sum_{n=0}^\infty P_n(z,q) t^n & = f(t) \\
&= \sum_{j=0}^\infty \frac{t^{2j} (1+t) q^{2j^2} (tzq;q^2)_j (tq/z;q^2)_j}
  { (t^2;q^2)_{2j+1}}\\
&= \sum_{j=0}^\infty t^{2j} (1+t) q^{2j^2} \sum_{h=0}^j (-tzq)^h q^{h^2-h}
\qbin{j}{h}{q^2} \sum_{i=0}^j (-tqz^{-1})^i q^{i^2-i} \qbin{i}{h}{q^2} \\
& \qquad \times \sum_{r=0}^\infty t^{2r} \qbin{2j+r}{r}{q^2}
\mbox{ \big(by~\eqref{qBT1} and~\eqref{qBT2}\big)}\\
&= \sum_{h,i,j,r\geq 0} t^{2j+h+i} (t^{2r} + t^{2r+1}) (-1)^{h+i}
q^{2j^2+h^2+i^2} z^{h-i}  \\
& \qquad \times \qbin{j}{h}{q^2}  \qbin{j}{i}{q^2}  \qbin{2j+r}{2j}{q^2} \\
&= \sum_{h,i,j,r\geq 0} t^{2j+h+i+s}  (-1)^{h+i}
q^{2j^2+h^2+i^2} z^{h-i}  \\
& \qquad \times \qbin{j}{h}{q^2}  \qbin{j}{i}{q^2}  \qbin{\lfloor \frac s 2 \rfloor+2j}{2j}{q^2}
\mbox{(where $s=2r$ or $s=2r+1$) }
\\
&=\sum_{n=0}^\infty t^n \sum_{h,i,j,k,\ell\geq 0}
(-1)^{h+i+\ell}
q^{2j^2+h^2+i^2+2k+2\ell} z^{h-i}   \qbin{j}{h}{q^2}  \qbin{i}{h}{q^2}
\\
& \qquad \times\qbin{j+ \lfloor \frac{n-h-i}{2} \rfloor}{2j}{q^2}\\
& \qquad\qquad\mbox{(where $n= 2j+h+i+s$)}.
\end{align*}
}
Compare co\"efficients of $t^n$ in the extremes to find
\begin{equation*} \label{R1finiteLHS}
 P_n(z,q) = \sum_{h,i,j\geq 0}
(-1)^{h+i} z^{h-i}
q^{2j^2+h^2+i^2}    \qbin{j}{h}{q^2}  \qbin{j}{i}{q^2}
\qbin{j+ \lfloor \frac{n-h-i}{2} \rfloor}{2j}{q^2}.
\end{equation*}

 Next, after some inspired guesswork, (see~\cite{AB90,S03,S04} for
details) we define the polynomials
\begin{multline*}
Q_n=Q_n(z,q) \\
:= \left\{
\begin{array}{ll} \sum_{k} z^{2k} q^{12k^2} \qbin{2m}{m+3k}{q^2} -
 z^{-2k-1} q^{12k^2+12k+3} \qbin{2m-1}{m+3k+1}{q^2}, &\mbox{  if $n=2m$}  \\
\sum_{k} z^{2k} q^{12k^2} \qbin{2m}{m+3k}{q^2} -
 z^{-2k-1} q^{12k^2+12k+3} \qbin{2m+1}{m+3k+2}{q^2}, &\mbox{  if $n=2m+1$}
 \end{array}
 \right. \\
  =\sum_{j=-\infty}^\infty (-1)^j z^j q^{3j^2} \qbin{n-1}{\lfloor \frac{n+3j-1}{2} \rfloor}{q^2}
 +\epsilon_n (z,q),
 \end{multline*}
 where
 {\allowdisplaybreaks
 \begin{equation*}
 \epsilon_n (z,q) \\= \left\{
 \begin{array}{ll}
\displaystyle{ \sum_{j=-\infty}^\infty z^{2j} q^{12j^2+6j+n} \qbin{n-1}{\frac{n+6j}{2}}{q^2}^*}  &\mbox{if $2\mid n$} \\
\displaystyle{-\sum_{j=-\infty}^\infty  z^{2j-1} q^{12j^2-6j+n} \qbin{n-1}{\frac{n+6j-3}{2}}{q^2}^*} &\mbox{if $2\nmid n$}
\end{array} \right. .
 \end{equation*}
 }
Our goal is to show that the $P_n(z,q)$ and $Q_n(z,q)$ are in fact one and the
same, thus giving us~\eqref{R1finite}.
We would like to use a computer implementation of the
$q$-Zeilberger algorithm~\cite{PWZ96,WZ90,WZ92,Z90,Z91} to simply show that the $Q_n$ satisfy
the recurrence~\eqref{R1rec}, and then upon checking that the
$Q_n$ satisfy the initial conditions~\eqref{R1ic}, we would be done.
Unfortunately, the implementations of the $q$-Zeilberger algorithm currently
available do not allow for direct input of summands as complex as those
under consideration here.   And the corresponding certificate function would likely
be rather horrendous.  Further, it is unlikely that the $q$-Zeilberger algorithm would
produce a minimal recurrence for the $Q_n$.  So, the traditional automated proof
would require a certain amount of pre-processing and post-processing.

\subsection{Identity~\ref{R1PartnerFinite} }
The derivation is analogous to that of Identity~\ref{R1finite}.
The recurrence and initial conditions are
\begin{multline}
  P_n = (1-q)P_{n-1} + (2q+q^n)P_{n-2} + \Big( q^2 - q -( zq^2 + z^{-1}q^3)q^{n-3} \Big)P_{n-3} \\+ (q^{n-1}-q^2) P_{n-4}
\end{multline}
with
\begin{multline*}
P_0 = 1,\quad P_1 = 1-q/z, \quad P_2 = 1+(1-z^{-1} )q+2q^2,\\ P_3 =
1+ (1-z^{-1}) q + (1-z-z^{-1} )q^2 -q^3/z - q^4/z .
\end{multline*}

\subsection{Identity~\ref{R2finite}}
The recurrence and initial conditions are
\begin{equation}
 P_n = (1+q-q^2+q^{2n-1})P_{n-1} + \Big( q^3+q^2 - q - (z+z^{-1} )q^{2n-2} \Big)
 P_{n-2} + (q^{2n-3}-q^3)P_{n-3}
\end{equation}
with
\[ P_0 = 1,\quad P_1 = 1+q, \quad P_2 = 1+q+(1-z-z^{-1})q^2 + q^4 .\]

\section{Challenge}
We leave it as a challenge to produce automated proofs for
Identities~\ref{R1finite}--~\ref{R2finite}.

\section*{Acknowledgments}
Many thanks to Doron Zeilberger for revolutionizing the way we approach the
discovery and proof of identities, especially those of the hypergeometric and
$q$-hypergeometric type.

.
\end{document}